\documentclass[a4paper,12pt]{article}
\usepackage[T2A]{fontenc} 
\usepackage[english]{babel}

\usepackage{mathrsfs}
\usepackage[top=2cm,left=3cm,right=1.5cm,bottom=2cm]{geometry}
\usepackage{amsmath}
\usepackage{amsfonts}
\usepackage{amssymb}
\usepackage{amsthm} 
\usepackage{enumerate}
\usepackage{hyperref}
\usepackage{alphalph}
\usepackage[subnum]{cases}

\numberwithin{equation}{subsection}

\theoremstyle{plain}
\newtheorem{theorem}{Theorem}[section]

\theoremstyle{definition}

\theoremstyle{remark}
\newtheorem{remark}[theorem]{Remark}

\numberwithin{equation}{section}

\begin{document}

%---------------------------------------------------------------------------

\begin{center}\textbf{About one inverse problem for the Sturm-Liouville operator}\end{center}
\vskip 0.3 cm
\begin{center}\textbf{B.\,N. Biyarov}\end{center}
\vskip 0.3 cm
%\begin{center}May 20, 2019\end{center}
\\

%%% ----------------------------------------------------------------------
%
\begin{abstract}
We consider the spectral problems for the  Sturm-Liouville operator generated by the Dirichlet, Neumann, Dirichlet-Neumann and Neumann-Dirichlet conditions. The necessary and sufficient condition for the coincidence of the spectrum of the Dirichlet-Neumann and Neumann-Dirichlet problems is proved. Also the necessary and sufficient condition for the coincidence of the spectrum, except zero,  of the Dirichlet and Neumann problems is proved. An application to periodic and anti-periodic problems is given.
\end{abstract}

\section{Introduction}
\label{sec:1} 
In the present paper we study the Sturm-Liouville operator
\[\widehat{L} = -\dfrac{d^2}{dx^2} + \ q(x), \]
in the Hilbert space $L^2(0, 1)$, where $\ q(x)$ is an arbitrary real-valued function of class $L^2(0, 1)$.
The closure in $L^2(0, 1)$ of the operator $\widehat{L}$ considered on $C^{\infty}[0, 1]$ is the maximal operator $\widehat{L}$ with the domain
\[D(\widehat{L}) = \{y \in L^2(0, 1): \ y,\ y' \in AC [0, 1], \  y'' - q(x)y \in L^2(0, 1)\}. \]
%-----------------------------------------
We consider the operator $L_D=\widehat{L}$ on the domain
\[D(L_D) = \{y \in D(\widehat{L}): \ y(0)=y(1)=0 \}, \]
the operator $L_N=\widehat{L}$ on the domain
\[D(L_N) = \{y \in D(\widehat{L}): \ y'(0)=y'(1)=0 \}, \]
the operator $L_{DN}=\widehat{L}$ on the domain
\[D(L_{DN}) = \{y \in D(\widehat{L}): \ y(0)=y'(1)=0 \}, \]
and the operator $L_{ND}=\widehat{L}$ on the domain
\[D(L_{ND}) = \{y \in D(\widehat{L}): \ y'(0)=y(1)=0 \}. \]
%------------------------------------------
Here we use subscripts $D$, $N$, $DN$ and $ND$ meaning Dirichlet, Neumann, Dirichlet-Neumann and Neumann-Dirichlet operators, respectively. By $\sigma(A)$ we denote the spectrum of the operator $A$.

One of the types of inverse problems for the Sturm-Liouville equation is to find some information about the potential $q(x)$ from a knowledge of the spectrum. Such problems have been studied by many authors (see \cite{Borg}, \cite{Hochstadt1}, \cite{Hochstadt2}, \cite{Hochstadt3}, \cite{Gesztesy1}, \cite{Gesztesy2}, \cite{Shkalikov},  \cite{Levitan}). In this paper, we want to find the properties of the potential $q (x)$ knowing that the spectrum of $L_{DN}$ and $L_{ND}$  coincides, as well as the spectrum of $L_{D}$ and $L_{N}$ coincides, except zero. It is known that these questions are directly related to the study of periodic or anti-periodic problems, the basicity of the system of root vectors, and also the questions of the coexistence of periodic solutions of Hill's equation.

The following theorems are the main results of this paper
%-------------Theorem1
\begin{theorem}
\label{theorem:1}
The spectrum of $L_{DN}$ coincides with the spectrum of $L_{ND}$ (i.e. $\sigma(L_{DN})=\sigma(L_{ND})$) if and only if $q(x)=q(1-x)$ on $[0, 1]$.
\end{theorem} 

%\setcounter{equation}{1}
%\def\theequation{\Alph{equation}}
%\def\theequation{\AlphAlph{\value{equation}}}
%\renewcommand{\theequation}{\Alph{equation}}

%-------------THEOREM2
\begin{theorem}\label{theorem:2}
The spectrum of $L_D$ coincides with the spectrum of  $L_N$, except zero (i.e. $\sigma(L_{D})\setminus \{0\}=\sigma(L_{N}) \setminus \{0\}$), and $0\in \sigma(L_N)$ if and only if 
\[
q_1(x)=\Bigg(\int \limits_{1}^{x}q_2(t)dt\Bigg)^2, \tag{BB}\label{MyEq}
\]
 where  $q_1(x)=(q(x)+q(1-x))/2$ and  $q_2(x)=(q(x)-q(1-x))/2$ on $[0, 1]$.
\end{theorem}
Theorem \ref{theorem:1} and Theorem \ref{theorem:2} will be proven in Section \ref{sec:3_ProofTh1} and \ref{sec:4_ProofTh2}, respectively.

%END Introduction
% ************************
%---Preliminaries
\section{Preliminaries}
\label{sec:2}
Consider the Sturm-Liouville equation in the Hilbert space $L^2(0, 1)$
\begin{equation}\label{eq:2.1}
    \widehat{L}y \equiv - y'' + q(x)y = \lambda^2y,
\end{equation}
where $q(x)$ is the real-valued function of class $L^2(0, 1)$.

%_---------------
By $c(x, \lambda)$ and $s(x, \lambda)$ we denote the fundamental system of solutions eqution \eqref{eq:2.1} corresponding to the initial conditions $c(0, \lambda) = s'(0, \lambda)=1$ and $c'(0, \lambda) = s(0, \lambda)=0$.
Then we have the representations (see \cite{Marchenco})

\begin{equation}\label{eq:2.2}
    \begin{cases}
        c(x, \lambda)=\cos{\lambda x} + \int\limits_{-x}^{x}K(x, t) \cos{\lambda t}dt, \\
        $$ s(x, \lambda)= \frac{\sin{\lambda x}}{\lambda} + \int\limits_{-x}^{x}K(x, t) \frac{\sin{\lambda t}}{\lambda}dt, $$ 
    \end{cases}
\end{equation}
in which $K(x, t) \in C(\Omega) \cap W_1^1(\Omega)$, where $\Omega=\{(x, t): 0 \leq x \leq 1, -x \leq t \leq x\}$, and $K(x, t)$ is the solution of the problem
\begin{equation}\label{eq:2.3}
    \begin{cases}
        K_{xx} - K_{tt} = q(x)K(x, t), \: \mbox{in} \  \Omega\\
        K(x, x) = \frac{1}{2}\int\limits_{0}^{x}q(t)dt, \: K(x, -x) = 0, \: x\in[0, 1]. 
    \end{cases}
\end{equation}
\\
We also consider the fundamental solutions of the eqution \eqref{eq:2.1} of the following form
\begin{equation}\label{eq:2.4}
    \begin{cases}
        y_1(x, \lambda)=\cos{\lambda (x-\frac{1}{2})} + \int\limits_{1-x}^{x}N(x, t) \cos{\lambda (t-\frac{1}{2})}dt, \\
        y_2(x, \lambda)= \frac {\sin{\lambda (x-\frac{1}{2})}}{\lambda} + \int\limits_{1-x}^{x}N(x, t) \frac{\sin{\lambda (t-\frac{1}{2})}}{\lambda}dt,  
    \end{cases}
\end{equation}
\\
in $\Omega_1 = \{ (x, t): \ 0 \leq x \leq 1, \ |t-\frac{1}{2}| \leq |x-\frac{1}{2}| \}$, with propoties $y_1(\frac{1}{2}, \lambda) = y'_{2}(\frac{1}{2}, \lambda)=1$, $y'_{1}(\frac{1}{2}, \lambda) = y_{2}(\frac{1}{2}, \lambda)=0$, where $N(x, t) = K(x-\frac{1}{2}, t-\frac{1}{2})$.
\\
Then
\begin{equation}\label{eq:2.5}
    \begin{cases}
        y_1(x, \lambda)= s'(\frac{1}{2}, \lambda)c(x, \lambda) - c'(\frac{1}{2}, \lambda)s(x, \lambda), \\
        y_2(x, \lambda)=  c(\frac{1}{2}, \lambda)s(x, \lambda) - s(\frac{1}{2}, \lambda)c(x, \lambda),  
    \end{cases}
\end{equation}
and

\begin{equation}\label{eq:2.6}
    \begin{cases}
        c(x, \lambda)= y'_{2}(0, \lambda)y_{1}(x, \lambda) - y'_{1}(0, \lambda)y_{2}(x, \lambda), \\
        s(x, \lambda)= -y_{2}(0, \lambda)y_{1}(x, \lambda) + y_{1}(0, \lambda)y_{2}(x, \lambda).  
    \end{cases}
\end{equation}
\\
From \eqref{eq:2.5} we have 
\[ y_{1}(0, \lambda) = s'(\frac{1}{2}, \lambda), \,\, y_{2}(0, \lambda) = -s(\frac{1}{2}, \lambda), \,\, y'_{1}(0, \lambda) = - c'(\frac{1}{2}, \lambda), \,\, y'_{2}(0, \lambda) = c(\frac{1}{2}, \lambda).\]
\\
Then the equation
\begin{equation}\label{eq:2.7}
    s'(1, \lambda) = c(1, \lambda), \: \forall \lambda \in \mathbb{C}
\end{equation}
and the equation
\begin{equation}\label{eq:2.8}
    y_{1}(0, \lambda) y'_{2}(1, \lambda) - y'_{1}(1, \lambda) y_{2}(0, \lambda) =  y'_{2}(0, \lambda) y_{1}(1, \lambda) - y'_{1}(0, \lambda) y_{2}(1, \lambda),
\end{equation}
will be equivalent.
\\
Then the equation
\begin{equation}\label{eq:2.9}
    c'(1, \lambda) = - \lambda^{2}s(1, \lambda),
\end{equation}
and the equation
\begin{equation}\label{eq:2.10}
    y'_{2}(0,\lambda) y_{1}^{'}(1, \lambda) - y'_{1}(0, \lambda) y'_{2}(1, \lambda) =  -\lambda^{2} [y_{1}(0, \lambda) y_{2}(1, \lambda) - y_{2}(0, \lambda) y_{1}(1, \lambda)],
\end{equation}
will be equivalent.
\\
We  calculate the following expression
\[ y_1(1-x, \lambda ) y_2(x, \lambda ) - y_1(x, \lambda ) y_2 (1 - x, \lambda ) \]
by virtue of \eqref{eq:2.4}
\[
    \begin{split}
        & y_1(1-x, \lambda)y_2(x, \lambda) - y_1(x, \lambda)y_2(1-x, \lambda) \\
        & = \bigg[
                \cos{\lambda(x-\frac{1}{2})}- \int\limits_{1-x}^{x} N(1-x, t)\cos\lambda(t-\frac{1}{2})dt 
            \bigg] \\
           & \cdot \bigg[
            \frac{\sin\lambda(x-\frac{1}{2})}{\lambda} + \int\limits_{1-x}^{x} N(x, t) \frac{\sin \lambda(t-\frac{1}{2})}{\lambda}dt
        \bigg] \\
        & +\bigg[
            \cos{\lambda(x-\frac{1}{2})} + \int \limits_{1-x}^{x}N(x, t) \cos{\lambda(t-\frac{1}{2})} dt
        \bigg] \\
        &\cdot  \bigg[
            \frac{\sin{\lambda(x-\frac{1}{2})}}{\lambda} + \int \limits_{1-x}^{x}N(1-x, t) \frac{\sin{\lambda(t-\frac{1}{2})}}{\lambda}dt
        \bigg]   \\
         & = \frac{1}{\lambda}\cos{\lambda(x-\frac{1}{2})} \sin{\lambda(x-\frac{1}{2})} +\frac{1}{\lambda} \int \limits_{1-x}^{x} N(x, t) \cos{\lambda (x-\frac{1}{2})} \sin{\lambda (t-\frac{1}{2})}dt 
    \end{split}
\]
\[
\begin{split}
        & - \frac{1}{\lambda} \int \limits_{1-x}^{x} N(1-x, t) \sin{\lambda (x-\frac{1}{2})} \cos{\lambda (t-\frac{1}{2})}dt \\
        & - \int \limits_{1-x}^{x} \int \limits_{1-x}^{x} N(x, t) N(1-x, \tau) \sin{\lambda(t-\frac{1}{2})} \cos{\lambda(\tau-\frac{1}{2})} dt d\tau  
    \end{split}
\]

\[
    \begin{split}
        & + \frac{1}{\lambda} \sin{\lambda (x-\frac{1}{2})} \cos{\lambda (x-\frac{1}{2})} + \int \limits_{1-x}^{x} N(1-x, t) \cos{\lambda(x-\frac{1}{2})} \sin{\lambda (t-\frac{1}{2}}) dt  \\
        & + \frac{1}{\lambda} \int \limits_{1-x}^{x} N(x, t) \sin{\lambda(x-\frac{1}{2})} \cos{\lambda(t-\frac{1}{2})} dt  \\
        & +\frac{1}{\lambda}\int \limits_{1-x}^{x} \int \limits_{1-x}^{x} N(1-x, \tau) N(x, t) \cos{\lambda(t-\frac{1}{2})} \sin{\lambda(\tau-\frac{1}{2})} dt d\tau 
    \end{split}
\]        
\[
    \begin{split}
        & =\frac{1}{\lambda}\sin{\lambda(2x-1)} + \frac{1}{\lambda} \int\limits_{1-x}^{x}N(x, t) \sin{\lambda(x+t-1)}dt   \\
        &  - \frac{1}{\lambda} \int\limits_{1-x}^{x}N(1-x, t)\sin{\lambda(x-t)}dt  -\frac{1}{\lambda} \int\limits_{1-x}^{x}\int\limits_{1-x}^{x}N(1-x, \tau)N(x,t) \sin{\lambda(t-\tau)}dt d \tau 
    \end{split}
\]        
\[
    \begin{split}        
         = \frac{1}{\lambda}&\sin{\lambda(2x-1)} +  \frac{1}{\lambda}\int \limits_{1-x}^{x} \big[ N(x, 1-t)-N(1-x, t) \big] \sin{\lambda (x-t)} dt \\
           & - \frac{1}{\lambda} \int \limits_{1-x}^{x} \int \limits_{1-x}^{x} N(1-x, \tau)N(x, t) \sin{\lambda (t-\tau)} dt d \tau  \\
        & = \frac{1}{\lambda} \sin{\lambda (2x-1)} + \frac{1}{\lambda} \int \limits_{0}^{2x-1} \big[ N(x, t-x+1) - N(1-x, x-t) \big] \sin{\lambda t} dt 
    \end{split}
\]        
\[
    \begin{split} 
        & - \frac{1}{\lambda} \int \limits_{1-x}^{x} dt \int \limits_{t-x}^{t+x-1} N(1-x, t- \tau) N(x, t) \sin{\lambda \tau} d \tau  \\
        &= \frac{1}{\lambda} \sin{\lambda (2x-1)} + \frac{1}{\lambda} \int \limits_{0}^{2x-1} \big[ N(x, t-x+1) - N(1-x, x-t) \big] \sin{\lambda t} dt 
    \end{split}
\]
\[
    \begin{split} 
        & - \int \limits_{0}^{2x-1} dt \int \limits_{t-(2x-1)}^{t} N(1-x, t-x+1- \tau) N(x, t-x+1) \sin{\lambda \tau} d \tau  \\
        & = \frac{1}{\lambda} \sin{\lambda (2x-1)} +\frac{1}{\lambda} \int \limits_{0}^{2x-1} \big[ N(x, t-x+1) - N(1-x, x-t) \big] \sin{\lambda t} dt 
    \end{split}
\]
\[
    \begin{split} 
        & - \frac{1}{\lambda} \int \limits_{-(2x-1)}^{0} d \tau \int \limits_{0}^{\tau + 2x - 1} N(1-x, t - x + 1 - \tau) N(x, t-x+1) \sin{\lambda \tau} dt  \\
        & - \frac{1}{\lambda} \int \limits_{0}^{2x-1} d \tau \int \limits_{\tau}^{2x-1} N(1-x, t-x+1-\tau) N(x, t-x+1) \sin{\lambda \tau} dt 
    \end{split}
\]

\[
    \begin{split} 
        & = \frac{1}{\lambda} \sin{\lambda (2x-1)} +\frac{1}{\lambda} \int \limits_{0}^{2x-1} \Big[ N(x, t-x+1) - N(1-x, x-t) \Big] \sin{\lambda t} dt  \\
        & + \frac{1}{\lambda}\int \limits_{0}^{2x-1}\sin{\lambda \tau} \bigg[ \int \limits_{0}^{2x-1 - \tau} N(1-x, \tau -x+t+1)N(x, t-x+1) dt  \\
        & - \int \limits_{\tau}^{2x-1} N(1-x, t-x+1-\tau)N(x, t-x+1) dt \bigg] d\tau.
    \end{split}
\]
Replace $\tau$ with $t$ in the last term of the equality. Then
\begin{equation}\label{eq:2.11}
    \begin{split}
    & y_1(1-x, \lambda)y_2(x, \lambda)-y_1(x, \lambda)y_2(1-x, \lambda)  \\
    & = \frac{1}{\lambda} \sin{\lambda (2x-1)} + \frac{1}{\lambda} \int \limits_{0}^{2x-1} \big[ N(x, t-x+1) - N(1-x, x-t) \big] \sin{\lambda t} dt  \\
    & + \frac{1}{\lambda} \int \limits_{0}^{2x-1} \sin{\lambda t} \bigg[ \int \limits_{0}^{2x-1-t} N(1-x, t-x+\tau+1)N(x, \tau-x+1) d\tau \\ 
    & - \int \limits_{t}^{2x-1} N(1-x, \tau-x+1-t) N(x, \tau-x+1) d \tau   \bigg] dt
    \end{split}
\end{equation}
We find the first derivative of the expression \eqref{eq:2.11}
\[
    \begin{split}
         \big[ y_1(1-x, &\lambda)y_2(x, \lambda) - y_1(x, \lambda)y_2(1-x, \lambda) \big]'  \\[6pt]
        & = 2\cos{\lambda(2x-1)} + \frac{2}{\lambda} \big[N(x, x) - N(1-x, 1-x) \big] \sin{\lambda(2x-1)}  \\
    \end{split}
\]
\[
    \begin{split}
        & +\frac{1}{\lambda} \int \limits_0^{2x-1} \Big[N'_1(x, t-x+1)-N'_2(x, t-x+1)+N'_1(1-x, x-t) \\
        & - N'_2(1-x, x-t)\Big]\sin{\lambda t}dt +  \frac{1}{\lambda} \int \limits_0^{2x-1} \sin{\lambda t} \bigg\{ 2N(1-x, x) N(x, x-t)  \\
        & +  \int \limits_0^{2x-1-t} \Big[ (-N'_1(1-x, t-x+\tau+1) - N'_2(1-x, t-x+\tau+1) ) 
    \end{split}
\]
\[
    \begin{split}
        & \cdot N(x, \tau-x+1) + N(1-x, t-x+\tau + 1) (N'_1(x, \tau-x+1)  \\[6pt]
        & - N'_2(x, \tau-x+1) ) \Big] d\tau - 2N(1-x, x-t)N(x, x) 
    \end{split}
\]
\[
    \begin{split}
        & - \int \limits_t^{2x-1} \Big[(-N'_1(1-x, \tau-x+1-t) - N'_2(1-x, \tau - x + 1 - t))  \\
        & \cdot N(x, \tau - x +1) + N(1-x, \tau - x + 1 - t) (N'_1(x, \tau - x + 1) \\[6pt]
        &- N'_2(x, \tau - x + 1) ) \Big] d\tau\bigg \}dt.
    \end{split}
\]
Since $N(1-x, x) = 0$ and denoting the expression in braces by $A(x, t)$, we have
\begin{equation}\label{eq:2.12}
    \begin{split}
    & \big[ y_1(1-x, \lambda)y_2(x, \lambda) - y_1(x, \lambda)y_2(1-x, \lambda) \big]'  \\[5pt]
   & = 2 \cos{\lambda(2x-1)} + \frac{1}{\lambda} \big[ N(x, x) - N(1-x, 1-x) \big] \sin{\lambda(x-1)}  \\
    & + \frac{1}{\lambda} \int \limits_0^{2x-1} \Big[ N'_1(x, t-x+1) - N'_2(x, t-x+1) + N'_1(1-x, x-t)  \\
    & - N'_2(1-x, x-t)\Big]\sin{\lambda t}dt + \frac{1}{\lambda} \int \limits_0^{2x-1} \sin{\lambda t} A(x, t)dt,
   \end{split}
\end{equation}
where
\[
   \begin{split}
    & A(x, t) = \int \limits_0^{2x-1-t} \Big\{ \big[ - N'_1(1-x, t-x+\tau+1) - N'_2(1-x, t-x+\tau+1) \big]  \\
    & \cdot N(x, \tau-x+1) + N(1-x, t-x+\tau+1) \cdot \big[N'_1(x, \tau-x+1)  \\[6pt]
     &- N'_2(x, \tau-x+1)\big] \Big\} d\tau - 2N(1-x, x-t)N(x, x)  \\
     & - \int \limits_t^{2x-1} \Big\{ \big[ - N'_1(1-x, \tau-x+1-t) - N'_2(1-x,  \tau-x+1-t) \big] \\
      & \cdot N(x, \tau - x + 1) +  N(1-x, \tau - x + 1 - t) \big[N'_1(x, \tau-x+1)  \\[6pt]
        & - N'_2(x, \tau-x+1) \big]\Big \} d\tau.
   \end{split}
\]
\\
Next, we transform the expression \eqref{eq:2.12}
\[
   \begin{split}
   & \big[ y_1(1-x, \lambda)y_2(x, \lambda) - y_1(x, \lambda)y_2(1-x, \lambda) \big]'  \\
   & = 2 \cos{\lambda(2x-1)} + 2 \big[ N(x, x) - N(1-x, 1-x) \big] \cdot \int \limits_0^{2x-1} \cos{\lambda t}dt  \\
     & + \int \limits_0^{2x-1} [ N'_1(x, t-x+1) - N'_2(x, t-x+1) + N'_1(1-x, x-t) \\
   & - N'_ 2(1-x, x-t) \cdot \int \limits_0^t \cos{\lambda \tau} d\tau + \int \limits_0^{2x-1} A(x, t) \cdot \int \limits_0^t \cos{\lambda \tau} d\tau dt 
   \end{split}
\]
\[
   \begin{split}
   & = 2\cos{\lambda(2x-1)} + 2 \big[ N(x, x) - N(1-x, 1-x) \big] \cdot \int \limits_0^{2x-1} \cos{\lambda t} dt  \\
   & +\int \limits_0^{2x-1} \cos{\lambda \tau} \int \limits_{\tau}^{2x-1} \Big[N'_1(x, t-x+1) - N'_2(x, t-x+1) + N'_1(1-x, x-t) 
   \end{split}
\]
\[
   \begin{split}
   & - N'_2(1-x, x-t) \Big] dt d\tau + \int \limits_0^{2x-1} \cos{\lambda \tau} \int \limits_{\tau}^{2x-1} A(x, t) dt d\tau  \\
    & = 2 \cos{\lambda(2x-1)} + 2 \big[ N(x, x) - N(1-x, 1-x) \big] \cdot \int \limits_{0}^{2x-1} \cos{\lambda t}dt 
   \end{split}
\]
\[
   \begin{split}
   & + \int \limits_0^{2x-1} \cos{\lambda t} \int \limits_t^{2x-1}  \Big[N'_1(x, \tau-x+1) - N'_2(x, \tau - x+1) + N'_1(1-x, x-\tau)  \\
   & - N'_2(1-x, x - \tau) \Big] d\tau dt + \int \limits_0^{2x-1} \cos{\lambda t}  \int \limits_t^{2x-1} A(x, \tau) d\tau dt.
   \end{split}
\]
Then we have
  \begin{equation}\label{eq:2.13}
 \big[ y_1(1-x, \lambda)y_2(x, \lambda) - y_1(x, \lambda)y_2(1-x, \lambda) \big]'  
 = 2 \cos{\lambda(2x-1)} +\int \limits_0^{2x-1} \cos{\lambda t} M(x, t)dt,
  \end{equation}
where

\[
   \begin{split}
    M(x, t) &= 2 \big[ N(x, x) - N(1-x, 1-x) \big] \\
   & +  \int \limits_t^{2x-1}  \Big[N'_1(x, \tau-x+1) - N'_2(x, \tau - x+1) + N'_1(1-x, x-\tau)  \\
   & - N'_2(1-x, x - \tau) \Big] d\tau + \int \limits_t^{2x-1} A(x, \tau) d\tau.
   \end{split}
\]
Throughout this paper, we used the notation $N'_1(x,t)$ and $N'_2(x, t)$, which denote the derivatives with respect to the first and the second arguments, respectively.
%###################################X
\section{Proof of Theorem \ref{theorem:1}}
\label{sec:3_ProofTh1}
From these equivalent equations \eqref{eq:2.7} and \eqref{eq:2.8} we have
\[
    \begin{split}
         s'(1, \lambda) - &c(1, \lambda) = y_1(0, \lambda)y'_2(1, \lambda) -   y'_1(1, \lambda)y_2(0, \lambda) \\[6pt]
        & +y'_1(0, \lambda)y_2(1, \lambda) -   y'_2(0, \lambda)y_1(1, \lambda) \\
           &= \int \limits_{0}^{1} (y_1(1-x, \lambda) y'_2(x, \lambda) + y'_1(1-x, \lambda) y_2(x, \lambda))'dx \\
        & =\int \limits_{0}^{1} (-y'_1(1-x, \lambda) y'_2(x, \lambda) + y_1(1-x, \lambda) y''_2(x, \lambda) \\
         & - y''_1(1-x, \lambda)y_2(x, \lambda) + y'_1(1-x, \lambda) y'_2(x, \lambda) ) dx  
    \end{split}
\]       
\[
    \begin{split}      
        & =\int \limits_{0}^{1} \big[ y_1(1-x, \lambda)y_2(x, \lambda)(q(x)-\lambda^2) - y_1(1-x, \lambda)y_2(x, \lambda)  (q(1-x)-\lambda^2)\big] dx  \\
       &= \int \limits_{0}^{1} (q(x) -q(1-x)) \cdot y_1(1-x, \lambda) y_2(x, \lambda)dx \\
       & = \int \limits_{0}^{\frac{1}{2}} (q(x) -q(1-x)) \cdot y_1(1-x, \lambda) y_2(x, \lambda)dx \\
    & + \int \limits_{\frac{1}{2}}^{1} (q(x) -q(1-x)) \cdot y_1(1-x, \lambda) y_2(x, \lambda)dx \\
    &= \int \limits_{\frac{1}{2}}^{1} (q(x) -q(1-x)) \cdot ( y_1(1-x, \lambda) y_2(x, \lambda) - y_1(x, \lambda) y_2(1-x, \lambda))dx.
    \end{split}
\]
Taking into account \eqref{eq:2.11} and making the denotation $g(x) = q(x)-q(1-x)$, we obtain
\[
    \begin{split}
     & \int \limits_{\frac{1}{2}}^{1} g(x) \Big\{ \frac{1}{\lambda}\sin{\lambda (2x-1)} + \frac{1}{\lambda} \int \limits_{0}^{2x-1} \big[N(x, t-x+1) - N(1-x, x-t)\big] \sin{\lambda t}dt \\
     &  + \frac{1}{\lambda} \int \limits_{0}^{2x-1} \sin{\lambda t} \Big[\int \limits_{0}^{2x-1-t} N(1-x, t-x+\tau+1) N(x, \tau-x+1) d\tau  \\
      & - \int \limits_{t}^{2x-1} N(1-x, \tau-x+1-t) N(x, \tau-x+1)d\tau \Big] dt\Big\}dx = 0.
     \end{split}
\]
%###################################XI
Then
\begin{equation}\label{eq:3.1}
    \begin{split}
        & \int \limits_{0}^{1} \frac{1}{\lambda} \sin{\lambda t} \Big\{ 2g(\frac{t+1}{2}) + \int \limits_{\frac{t+1}{2}}^{1} g(x) \Big[N(x, t-x+1)  \\
        & - N(1-x, x-t) + \int \limits_{0}^{2x-1-t} N(1-x, t-x+\tau+1) N(x, \tau - x+1) d\tau  \\
        & - \int \limits_{t}^{2x-1} N(1-x, \tau-x+1-t) N(x, \tau - x+1) d\tau \Big] dx \Big\}dt=0.
     \end{split}
\end{equation}
If $g(x) = q(x) - q(1-x) = 0$ in \eqref{eq:3.1}, then $s'(1, \lambda) = c(1, \lambda)$ for all $\lambda \in \mathbb{C}$. Hence $\sigma(L_{DN}) = \sigma(L_{ND})$.
We now prove the converse of this theorem.
Let $\sigma(L_{DN}) = \sigma(L_{ND})$  and $\{\lambda_n\}_1^{\infty} = \{ \lambda \in \mathbb{C}: \  c(1, \lambda) = s'(1, \lambda) = 0\}$.

The system $\{ \sin{\lambda_{n}t} \} _1^{\infty}$ is complete in $L^2(0, 1)$ because it is equivalent \cite[p.10]{Marchenco} to the complete system $\{ s(t, \lambda_{n}) \}_1^{\infty}$ of eigenfunctions of $L_{DN}$.
Therefore, from \eqref{eq:3.1}  we get
\[
    \begin{split}
        & 2g(\frac{t+1}{2}) + \int \limits_{\frac{t+1}{2}}^{1} g(x)\Big[N(x, t-x+1) - N(1-x, x-t)  \\
        &+ \int \limits_{0}^{2x-1-t} N(1-x, t-x+1+\tau) N(x, \tau - x + 1) d\tau  \\
        & -\int \limits_{t}^{2x-1} N(1-x, \tau-x+1-t) N(x, \tau - x + 1) d\tau \Big] dx = 0. 
    \end{split}
\]
Making the change of variable $x = \xi$ in this integral, we obtain
\[
    \begin{split}
        & 2g(\frac{t+1}{2}) + \int \limits_{\frac{t+1}{2}}^{1} g(\xi)\Big[N(\xi, t-\xi+1) - N(1-\xi, \xi-t)  \\
        & +\int \limits_{0}^{2\xi-1-t} N(1-\xi, t-\xi+\tau+1) N(\xi, \tau - \xi + 1) d\tau \\
        & - \int \limits_t^{2 \xi - 1} N(1- \xi, \tau - \xi +1 - t) N(\xi, \tau - \xi +1) d\tau\Big] d\xi = 0. 
    \end{split}
\]
Denoting $x = \frac{t+1}{2}$,  we obtain
\begin{equation}\label{eq:3.2}
    \begin{split}
       & g(x) + \int \limits_x^{1} g(\xi) \Big[ N(\xi, 2x - \xi) - N(1 -\xi, \xi - 2x +1)  \\
       & + \int \limits_0^{2(\xi - x)} N(1- \xi, 2x - \xi + \tau) N(\xi, \tau - \xi +1) d\tau   \\
       &- \int \limits_{2x - 1}^{2\xi - 1} N(1- \xi, \tau- \xi - 2x + 2) N(\xi, \tau - \xi +1) d \tau \Big] d\xi = 0.
    \end{split}
\end{equation}
Equation \eqref{eq:3.2} has only the zero solution since some power of the nonlinear operator $(Ug)(x) = \int \limits_{x}^{1} g(\xi) G(x, \xi) d\xi$, is a contracting operator due to the fact that $|G(x, \xi)| \leq c$ by virtue of the estimate (see \cite[p.28]{Marchenco}) $|N(x, t)| \leq c$, where 
\[
    \begin{split}
        & G(x, \xi) = N(\xi, 2x - \xi) - N(1-\xi, \xi - 2x + 1)  \\
        & + \int \limits_{0}^{2(\xi -x)} N(1-\xi, 2x - \xi + \tau) N(\xi, \tau - \xi + 1) d\tau  \\
        & - \int \limits_{2x-1}^{2\xi -1)} N(1-\xi, \tau- \xi -2x + 2 ) N(\xi, \tau - \xi + 1)d \tau .
    \end{split}
\]
Thus, we have proved the Theorem \ref{theorem:1}.
\section{Proof of Theorem \ref{theorem:2}}
\label{sec:4_ProofTh2}
From these equivalent equations \eqref{eq:2.9} and \eqref{eq:2.10} we have.
\[
    \begin{split}
    c'(1, \lambda)& + \lambda^2 s(1, \lambda) = y'_2(0, \lambda)y'_1(1, \lambda) - y'_1(0, \lambda)y'_2(1, \lambda)  \\[7pt]
    & +\lambda^2 y_1(0, \lambda) y_2(1, \lambda) - \lambda^2 y_2(0, \lambda) y_1(1, \lambda) \\
     & = \int \limits_{0}^{1} \big[ y'_1(x, \lambda)y'_2(1-x, \lambda) - \lambda^2y_1(x, \lambda)y_2(1-x, \lambda) \big]' dx  \\
    & = \int \limits_{0}^{2} \Big[ y''_1(x, \lambda)y'_2(1-x, \lambda) - y'_1(x, \lambda)y''_2(1-x, \lambda) \\
     & - \lambda^2y'_1(x, \lambda)y_2(1-x, \lambda) + \lambda^2y_1(x, \lambda)y'_2(1-x, \lambda)\Big ]dx   \\
    & = \int \limits_{0}^{1}\Big [ y_1(x, \lambda) y'_2(1-x, \lambda) (q(x) - \lambda^2) - y'_1(x, \lambda) y_2(1-x, \lambda) (q(1-x) - \lambda^2) \\
     & - \lambda^2 y'_1(x, \lambda) y_2(1-x, \lambda) + \lambda^2 y_1(x, \lambda) y'_2(1-x, \lambda)\Big ] dx  \\
    & =\int \limits_{0}^{1} \big[ y_1(x, \lambda)y'_2(1-x, \lambda)q(x) - y'_1(x, \lambda)y_2(1-x, \lambda)q(1-x) \big] dx = 0.
    \end{split}
\]
We denote $q(x)=q_1(x) + q_2(x)$, where $q_1(x)=q_1(1-x)$, $q_2(x)=-q_2(1-x)$, on $[0, 1]$.
\[
    \begin{split}
     c'(1, \lambda) + \lambda^2s(1, \lambda)& = \int \limits_0^1 q_1(x) \big[y_1(x, \lambda)y'_2(1-x, \lambda) - y'_1(x, \lambda)y_2(1-x, \lambda)\big]dx  \\
    & +\int \limits_0^1 q_2(x) \big[y_1(x, \lambda)y'_2(1-x, \lambda) + y'_1(x, \lambda)y_2(1-x, \lambda)\big]dx  
    \end{split}
\]
\[
\begin{split}
  &= \int \limits_{0}^{\frac{1}{2}} q_1(x) \big[y_1(x, \lambda)y'_2(1-x, \lambda) - y'_1(x, \lambda)y_2(1-x, \lambda)\big]dx  \\
     & +  \int \limits_{\frac{1}{2}}^1 q_1(x) \big[y_1(x, \lambda)y'_2(1-x, \lambda) - y'_1(x, \lambda)y_2(1-x, \lambda)\big]dx  \\
     & +  \int \limits_{0}^{\frac{1}{2}} q_2(x)\big [y_1(x, \lambda)y'_2(1-x, \lambda) + y'_1(x, \lambda)y_2(1-x, \lambda)\big]dx \\
       & + \int \limits_{\frac{1}{2}}^{1} q_2(x) \big[ y_1(x, \lambda ) y'_2(1-x, \lambda) + y'_1(x, \lambda ) y_2(1-x, \lambda)\big] dx 
  \end{split}
\]
\[
    \begin{split}
      & = \int \limits_{\frac{1}{2}}^{1} q_1(x) \Big[ y_1(x, \lambda ) y'_2(1-x, \lambda) - y'_1(x, \lambda ) y_2(1-x, \lambda)  \\
        & + y_1(1-x, \lambda ) y'_2(x, \lambda) - y'_1(1-x, \lambda ) y_2(x, \lambda)\Big] dx \\
        & + \int \limits_{\frac{1}{2}}^{1} q_2(x) \Big[ y_1(x, \lambda ) y'_2(1-x, \lambda) + y'_1(x, \lambda ) y_2(1-x, \lambda)  \\
        & - y_1(1-x, \lambda ) y'_2(x, \lambda) - y'_1(1-x, \lambda ) y_2(x, \lambda)\Big] dx 
    \end{split}
\]
\[
    \begin{split}
      & = \int \limits_{\frac{1}{2}}^{1} q_1(x) \big[ y_1(1-x, \lambda ) y_2(x, \lambda) - y_1(x, \lambda ) y_2(1-x, \lambda)\big]'dx  \\
        & + \int \limits_{1}^{x} q_2(t)dt \Big[ y_1(x, \lambda ) y'_2(1-x, \lambda) + y'_1(x, \lambda ) y_2(1-x, \lambda) \\
        & - y_1(1-x, \lambda) y'_2(x, \lambda) - y'_1(1-x, \lambda) y_2(x, \lambda) \Big] \Big|_{\frac{1}{2}}^1 
         \\
        & - \int \limits_{\frac{1}{2}}^{1} \int \limits_{1}^{x} q_2(t)dt \Big[ y'_1(x, \lambda ) y'_2(1-x, \lambda)-y_1(x, \lambda) y''_2(1-x, \lambda) \\
          & +y''_1(x, \lambda) y_2(1-x, \lambda) - y'_1(x, \lambda) y'_2(1-x, \lambda) + y_1'(1-x, \lambda)y'_2(x, \lambda)  \\[8pt]
        & - y_1(1-x, \lambda)y''_2(x, \lambda) + y''_1(1-x, \lambda)y_2(x, \lambda) - y'_1(1-x, \lambda)y_2(x, \lambda)\Big ] dx 
    \end{split}
\]
\[
    \begin{split}
        & = \int \limits_{\frac{1}{2}}^{1} q_1(x) \big[ y_1(1-x, \lambda) y_2(x, \lambda) - y_1(x, \lambda)y_2(1-x, \lambda)\big]' dx  \\
        &-\int \limits_{\frac{1}{2}}^{1} \int \limits_{1}^{x} q_2(t)dt \Big[-y_1(x, \lambda) y_2(1-x, \lambda) (q(1-x) - \lambda^2) \\
          & +y_1(x, \lambda) y_2(1-x, \lambda) (q(x) - \lambda^2) - y_1(1-x, \lambda)y_2(x, \lambda) (q(x) - \lambda^2)  \\[8pt]
   & + y_1(1-x, \lambda)y_2(x, \lambda)(q(1-x) - \lambda^2) \Big]dx     
    \end{split}
\]
\[
    \begin{split}
    & =  \int \limits_{\frac{1}{2}}^{1} q_1(x) \big[ y_1(1-x, \lambda) y_2(x, \lambda) - y_1(x, \lambda)y_2(1-x, \lambda)\big]' dx \\
        & +\int \limits_{\frac{1}{2}}^{1} q_2(x) ( q(x) - q(1-x) )\big [ y_1(1-x, \lambda)y_2(x, \lambda) - y_1(x, \lambda)y_2( 1-x, \lambda )\big ] dx 
    \end{split}
\]
\[
    \begin{split}
        & = \int \limits_{\frac{1}{2}}^{1} q_1(x)\big[ y_1(1-x, \lambda)y_2(x, \lambda) - y_1(x, \lambda)y_2(1-x, \lambda) \big]' dx  \\
         &+ \int \limits_{\frac{1}{2}}^{1} \int \limits_{1}^{x} q_2(t)dt 2q_2(x)\big[ y_1(1-x, \lambda)y_2(x, \lambda) - y_1(x, \lambda)y_2(1-x, \lambda)\big] dx  \\
    \end{split}
\]
\[
    \begin{split}
        & = \int \limits_{\frac{1}{2}}^{1} q_1(x)\big[ y_1(1-x, \lambda)y_2(x, \lambda) - y_1(x, \lambda)y_2(1-x, \lambda)\big ]' dx  \\
         & - \int \limits_{\frac{1}{2}}^{1} 2 \int \limits_{1}^{x} q_2(t) \int \limits_{1}^{t} q_2(\tau) d\tau dt\big [ y_1(1-x, \lambda)y_2(x, \lambda) - y_1(x, \lambda)y_2(1-x, \lambda)\big]' dx  \\
         &= \int \limits_{\frac{1}{2}}^{1} Q(x)\big [ y_1(1-x, \lambda)y_2(x, \lambda) - y_1(x, \lambda)y_2(1-x, \lambda) \big]' dx = 0,
    \end{split}
\]
where $Q(x) = q_1(x) - 2 \int \limits_{1}^{x} q_2(t) \int \limits_{1}^{t} q_2(\tau)d\tau dt$=$q_1(x)-\Bigg(\int \limits_{1}^{x}q_2(t)dt\Bigg)^2$.
By virtue of \eqref{eq:2.13}, we have
\[\int \limits_{\frac{1}{2}}^{1} Q(x) \Big[2 \cos{\lambda (2x-1)} + \int \limits_{0}^{2x-1} \cos{\lambda t}M(x, t)dt \Big]dx = 0.\]
Hence
  \begin{equation}\label{eq:4.1} \int \limits_{0}^{1} \cos{\lambda t} \Big[ 2Q(\frac{t+1}{2}) +\int \limits_{\frac{t+1}{2}}^{1} Q(x)M(x, t)dx \Big] dt = 0.
  \end{equation}
If $Q(x)=0$ in \eqref{eq:4.1}, then $c'(1, \lambda) = - \lambda^{2}s(1, \lambda)$ for all $\lambda \in \mathbb{C}$.
Hence the equality $\sigma(L_{D})\setminus \{0\}=\sigma(L_{N}) \setminus \{0\}$ and $0\in \sigma(L_N)$ are obvious. 

Let $\sigma(L_{D})\setminus \{0\}=\sigma(L_{N}) \setminus \{0\}$, i.e.,  $\{ \lambda_{n} \}_1^{\infty} = \{ \lambda \in \mathbb{C}\setminus \{0\} : c'(1, \lambda) = - \lambda^2 s(1, \lambda) = 0 \}$
and $\lambda_0=0\in \sigma(L_N)$.
The system $\{ \cos{\lambda_{n}t} \} _0^{\infty}$ is complete in $L^2(0, 1)$ because it is equivalent \cite[p.10]{Marchenco} to the complete system $\{ c(t, \lambda_{n}) \}_0^{\infty}$ of eigenfunctions of $L_N$.
Therefore, from \eqref{eq:4.1}  we get
\[
    Q(\frac{t+1}{2}) + \frac{1}{2} \int \limits_{\frac{t+1}{2}}^1 Q(x)M(x, t)dx = 0.
\]
Making the change of variable $x= \xi$ in this integral, we obtain
\[
    Q(\frac{t+1}{2}) + \frac{1}{2} \int \limits_{\frac{t+1}{2}}^1 Q(\xi)M(\xi, t)d\xi = 0.
\]
If we  denote $x=\frac{t+1}{2}$, then we have
\begin{equation}\label{eq:4.2}
     Q(x)+\frac{1}{2}\int \limits_{{x}}^1 Q(\xi)M(\xi, 2x-1)d\xi = 0.
\end{equation}
Equation \eqref{eq:4.2} has only the zero solution since some power of the nonlinear operator
\[
    (UQ)(x) = \frac{1}{2}\int \limits_{x}^1 Q(\xi)M(\xi, 2x-1)d\xi
\]
is a contracting operator the to fact that $|M(\xi, 2x-1)| \leq c$ by virtue of the estimate (see \cite[p.28]{Marchenco}). Thus, we have proved the Theorem \ref{theorem:2}.
%####################
%
\section{Applications}
\label{sec:5}
In this section we consider the operator $L_P = \widehat{L}$ on the domain
\[
D(L_P)=\{ y \in D(\widehat{L}): y(0)=y(1), \, y'(0)=y'(1) \}
\]
and the operator $L_{AP}=\widehat{L}$ on the domain
\[
D(L_{AP})=\{ y \in D(\widehat{L}): y(0)=-y(1), \, y'(0)=-y'(1) \}.
\]
Let  $q(x)=q(1-x)$ on $[0, 1]$. It follows from \eqref{eq:2.5} and \eqref{eq:2.6} that
\begin{enumerate}
    \item[(a) ] $s'(1, \lambda) = c(1, \lambda)$,
    \item[(b) ] $s(1, \lambda)=2s(\frac{1}{2}, \lambda)s'(\frac{1}{2}, \lambda)$,
    \item[(c) ] $c'(1, \lambda)=2c(\frac{1}{2}, \lambda)c'(\frac{1}{2}, \lambda)$,
    \item[(d) ] $c(1, \lambda)=c(\frac{1}{2}, \lambda)s'(\frac{1}{2}, \lambda) + s(\frac{1}{2}, \lambda)c'(\frac{1}{2}, \lambda) = 1 + 2s(\frac{1}{2}, \lambda)c'(\frac{1}{2}, \lambda) = 2c(\frac{1}{2}, \lambda)s'(\frac{1}{2}, \lambda) - 1$,
    \item[(e) ] $s'(1, \lambda)=s(\frac{1}{2}, \lambda)c'(\frac{1}{2}, \lambda)+s'(\frac{1}{2}, \lambda)c(\frac{1}{2}, \lambda) = 1 + 2s(\frac{1}{2})c'(\frac{1}{2}, \lambda) = 2c(\frac{1}{2}, \lambda)s'(\frac{1}{2}, \lambda)-1$,
\end{enumerate}
for all $\lambda \in \mathbb{C}$. Then we obtain 
\[
\sigma(L_{P}) = \{ \lambda \in \mathbb{C}: s(\frac{1}{2}, \lambda) c'(\frac{1}{2}, \lambda) = 0 \},
\]
\[
\sigma(L_{AP}) = \{ \lambda \in \mathbb{C}: c(\frac{1}{2}, \lambda) s'(\frac{1}{2}, \lambda) = 0 \}.
\]
For the eigenvalues of  $L_P$ there are following three possible cases

\begin{enumerate}
    \item[(i) ] $s(\frac{1}{2}, \lambda) = 0, \ c'(\frac{1}{2}, \lambda) \neq 0$, we denote them by $\{ \lambda_n^s \}_1^\infty$,
    \item[(ii) ] $c'(\frac{1}{2}, \lambda) = 0, s(\frac{1}{2}, \lambda) \neq 0$, we denote them by $\{ \lambda_n^{c'} \}_{0}^{\infty}$,
    \item[(iii) ] $s(\frac{1}{2}, \lambda) = c'(\frac{1}{2}, \lambda) = 0$, we denote them by $\{ \lambda_n^{sc'} \}_{1}^{\infty}$.
\end{enumerate}

The following eigenvectors correspond to each of the indicated cases, respectively
\begin{enumerate}
    \item{
        $y_2(x, \lambda_n^s) = c(\frac{1}{2}, \lambda_n^s)s(x, \lambda_n^s), $
    }
    \item{
         $y_1(x, \lambda_n^{c'}) = s'(\frac{1}{2}, \lambda_n^{c'})c(x, \lambda_n^{c'}), $
    }
    \end{enumerate}
\begin{numcases}{3.} 
            y_1(x, \lambda_n^{sc'}) = s'(\frac{1}{2}, \lambda_n^{sc'})c(x, \lambda_n^{sc'}) \label{eq:5.1a} \\ 
            y_2(x, \lambda_n^{sc'}) = c(\frac{1}{2}, \lambda_n^{sc'})s(x, \lambda_n^{sc'})  \label{eq:5.1b}
 \end{numcases}  
We consider the operator $L_{D(\frac{1}{2})}=\widehat{L}$ on $[0,\frac{1}{2}]$ on the domain
\[
    D(L_{D(\frac{1}{2})}) = \{ y \in D(\widehat{L}): y(0)=y(\frac{1}{2})=0 \}
\]
and  the operator $L_{N(\frac{1}{2})}=\widehat{L}$ on $[0, \frac{1}{2}]$ on the domain
\[
    D(L_{N(\frac{1}{2})})=\{ y \in  D(\widehat{L}): y'(0)=y'(\frac{1}{2})=0 \}. 
\]
The condition \eqref{MyEq}  on $\big[0, \frac{1}{2} \big]$ we rewrite in the form
\[
q_1(x)=\Bigg(\int \limits_{\frac{1}{2}}^{x}q_2(t)dt\Bigg)^2, \tag{B}\label{MyEq2}
\]
where $q_1(x)=(q(x)+q(\frac{1}{2}-x))/2$ and  $q_2(x)=(q(x)-q(\frac{1}{2}-x))/2$ on $\big[ 0, \frac{1}{2} \big]$.
By Theorem \ref{theorem:2}, the condition \eqref{MyEq2} is necessary and sufficient for the coincidence of the spectrum of  $L_{D(\frac{1}{2})}$ and $L_{N(\frac{1}{2})}$, except zero, (i.e. $\sigma(L_{D(\frac{1}{2})})\setminus \{0\}=\sigma(L_{N(\frac{1}{2})}) \setminus \{0\}$), and $0\in \sigma(L_{N(\frac{1}{2}})$.

\begin{theorem}\label{theorem:5.1}
Let $q(x)=q(1-x)$ on $[0, 1]$. Then, the whole spectrum of the operator $L_P$ , except the lowest, consists only of eigenvalues with multiplicity two if and only if  the condition \eqref{MyEq2} holds. Moreover, one of the eigenfunctions  \eqref{eq:5.1a}, corresponding to the eigenvalue $\lambda_n^{sc'}$, is even on $[ 0, 1]$  and satisfies the condition of the Neuman problem on $[0, 1]$, and the other \eqref{eq:5.1b} is odd on $\big[ 0, 1 \big]$ and satisfies the condition of the Dirichlet problem on $[ 0, 1]$.
\end{theorem}

For the eqgenvalues of the operator $L_{AP}$ there are following three possible cases:
\begin{enumerate}[(i)]
    \item $c(\frac{1}{2}, \lambda)=0, s'(\frac{1}{2}, \lambda)\neq 0 $, we denote them by $ \{ \lambda_n^c \}_1^\infty$,
    \item $s'(\frac{1}{2}, \lambda) = 0$, $c(\frac{1}{2}, \lambda) \neq 0$, we denote them by $ \{ \lambda_n^{s'} \}_1^\infty$,
    \item $c(\frac{1}{2}, \lambda) = s'(\frac{1}{2}, \lambda) = 0$, we denote them by $ \{ \lambda_n^{cs'} \}_1^\infty$.
\end{enumerate}
The following eigenfunctions correspond to each of the indicated cases, respectively
\begin{enumerate}
    \item $y_2(x, \lambda_n^{c}) = -s(\frac{1}{2}, \lambda_n^c) c(x, \lambda_n^c)$,
    \item $y_1(x, \lambda_n^{s'}) = -c'(\frac{1}{2}, \lambda_n^{s'}) s(x, \lambda_n^{s'})$,
\end{enumerate}
\begin{numcases}{3.} 
          y_1(x, \lambda_n^{cs'}) = -c'(\frac{1}{2}, \lambda_n^{cs'}) s(x, \lambda_n^{cs'}),  \label{eq:5.2a} \\
          y_2(x, \lambda_n^{cs'}) = -s(\frac{1}{2}, \lambda_n^{cs'}) c(x, \lambda_n^{cs'}). \label{eq:5.2b} 
\end{numcases}

We consider the operator $L_{DN(\frac{1}{2})}=\widehat{L}$ on $[0,\frac{1}{2}]$ on the domain
\[
    D(L_{DN(\frac{1}{2})}) = \{ y \in D(\widehat{L}): y(0)=0, \, y'(\frac{1}{2})=0 \}
\]
and  the operator $L_{ND(\frac{1}{2})}=\widehat{L}$ on $[0, \frac{1}{2}]$ on the domain
\[
    D(L_{ND(\frac{1}{2})})=\{ y \in  D(\widehat{L}): y'(0)=0, \, y(\frac{1}{2})=0 \}.
\]
By  Theorem \ref{theorem:1}, the condition $q(x)=q(\frac{1}{2}-x)$ on $[0, \frac{1}{2}]$ is necessary and sufficient for the coincidence of the spectrum of  $L_{DN(\frac{1}{2})}$ and $L_{ND(\frac{1}{2})}$ (i.e. $\sigma(L_{DN(\frac{1}{2})} = L_{ND(\frac{1}{2})}$).
Thus, we have proved the following
\begin{theorem}\label{theorem:5.2} 
    Let $q(x)=q(1-x)$ on $[0, 1]$. Then, the whole spectrum of $L_{AP}$  consists only of  eigenvalues with multiplicaty two if and only if the  condition $q(x)=q(\frac{1}{2}-x)$ on $[0, \frac{1}{2}]$ holds. Moreover, one of the eigenfunctions \eqref{eq:5.2a} corresponding to the eigenvalue $\lambda_n^{cs'}$, is even on $[0, 1]$ and satisfies the condition Dirichlet on $[0, 1]$, and the other \eqref{eq:5.2b}  is odd on $[0, 1]$ and satisfies the condition of the Neuman problem on $[0, 1]$.
\end{theorem}

\begin{remark}\label{rem5.3}
Theorems \ref{theorem:1}, \ref{theorem:2},  \ref{theorem:5.1}, and  \ref{theorem:5.2} remain also valid for complex-valued $q(x)$ in $L^2(0,1)$.
\end{remark}

Indeed, in the proof of these theorems, we have not used the reality of $q(x)$.

\begin{remark}\label{rem5.4}
It is known that in the particular case $q(x)=q=const$, we have
\[\sigma(L_{D})=\sigma(L_{N})\setminus \{q\}.\]
\end{remark}

\begin{flushleft}
   Bazarkan Nuroldinovich Biyarov \\
   Faculty of Mechanics and Mathematics\\
   L.\,N.\,Gumilyov Eurasian National University \\
   Satpayev Str., 2\\
   010008 Nur-Sultan, Kazakhstan\\
   E-mail: bbiyarov@gmail.com
\end{flushleft} 

\end{document}